\documentclass[a4paper,11pt]{article}

\usepackage{color}
\usepackage{graphicx}
\usepackage{amsmath}
\usepackage{amsfonts}
\usepackage{amssymb}
\usepackage{amsthm}
\usepackage{epstopdf}
\usepackage{setspace}
\usepackage{fullpage}
\usepackage{url}

\theoremstyle{remark}
\newtheorem{remark}{Remark}
\newtheorem{algorithm}{Algorithm}
\newtheorem{problem}{Problem}
\newtheorem{lemma}{Lemma}
\newtheorem{corollary}{Corollary}

\newcommand{\intO}{\int_\Omega}

\DeclareMathOperator*{\argmin}{argmin}

\title{Verification of functional a posteriori \\ error estimates for obstacle problem in 1D}

\author{P. HARASIM, J. VALDMAN
}

\begin{document}
\maketitle

\begin{abstract}
We verify functional a posteriori error estimate for obstacle problem proposed by Repin.
Simplification into 1D allows for the construction of a nonlinear benchmark for which an exact solution of the obstacle problem can be derived. Quality of a numerical approximation obtained by the finite element method is compared with the exact solution and the error of approximation is bounded from above by a majorant error estimate. The sharpness of the majorant error estimate is discussed.
\end{abstract}

\section{Introduction}
Obstacle problems are one of the key problems in continuum mechanics. Their mathematical models based on variational inequalities are well established (we refer to classical works \cite{GLT, HHNL,KO}). Numerical treatment of a obstacle problem is obtained by the finite element method and a solution of a quadratic minimization problems with constrains. It was traditionally tackled by the Uzawa method, the interior point method, the active set method with gradient splitting and the semi-smooth Newton method among others \cite{DO, Ul}. \\

A priori analysis providing asymptotic estimates of the quality of finite elements approximations converging toward the exact solution was studied for obstacle problems e.g. in \cite{BHR, FALK}. For the survey of the most important techniques in a posteriori analysis (such as residual, gradient averaging or equilibration methods) we refer to the monographs \cite{AiOd, BabSt, BanRa}. Particular a posteriori estimates for variational inequalities including a obstacle problem are reported e.g. in \cite{BHS, CaMe, ZVKG} among others. \\

Our goal is to verify guaranteed functional a posteriori estimates expressed in terms of functional majorants derived by Repin \cite{NeRep, ReGruyter}. The functional majorant upper bounds are essentially different with respect to known a posteriori error estimates mentioned above. The estimates are obtained with the help of variational (duality) method which was developed in \cite{Re1,Re2} for convex variational problems. The method was applied to various nonlinear models including those associated with variational inequalities \cite{Re4}, in particular problems with obstacles \cite{BuRe}, problems generated by plasticity theory \cite{FuRe, ReVa2} and problems with nonlinear boundary conditions \cite{ReVa}.\\

The obstacle problem is formulated and analyzed in two dimensions, however numerical experiments are considered in one dimension only. Then, we are easily able to construct an analytical benchmark with an exact solution of the nonlinear obstacle problem and evaluate integrals in numerical tests exactly. \\

Outline of the paper is as follows. In Section 2, we formulate a constrained minimization problem and introduce a perturbed minimization problem including its basic properties. A derivation and further analysis of error estimates in term of a functional majorant is explained in Section 3. A method of majorant minimization is also included there.
A benchmark with known analytical solution is discussed in Section 4. Numerical tests performed in Matlab are reported in Section 5.

\section{Formulation of obstacle problem}
Let $\Omega\subset\mathbb{R}^2$
is a bounded domain with Lipschitz continuous boundary $\partial\Omega$. Let
$V$ stands for the standard Sobolev space $H^{1}(\Omega)$ and $V_{0}$ denote its subspace $H^{1}_{0}(\Omega)$, consisting
of functions whose trace on $\partial\Omega$ is zero.
We consider the obstacle problem, described by the following minimization problem:
\begin{problem}[Minimization problem] \label{prob1} Find $u \in K$ satisfying
\begin{equation}
J(u)=\inf_{v \in K} J(v),   \label{2rov1}
\end{equation}
where the energy functional reads
\begin{equation}
J(v) := \frac{1}{2}\intO\nabla v\cdot\nabla v \,\textrm{d}x-\intO fv\,\textrm{d}x      \label{2rov2}
\end{equation}
and the admissible set is defined as
$$
K:=\bigl\{v \in V_{0}:\,
v(x) \geq \phi(x) \mbox{ a.e.}\;\textrm{in}\; \Omega \},
$$
where $f\in L^{2}(\Omega)$ and $\phi\in V$ such that $\phi\not\in V_{0}$ and $\phi(x)<0 \mbox{ a.e.}\;\textrm{in}\; (\Omega)$.
\end{problem}

Problem \ref{prob1} is a quadratic minimization problem with a convex constrain and the existence of its minimizer is guaranteed by the Lions-
Stampacchia Theorem \cite{LiSt}. It is equivalent to the following variational inequality: Find $u \in K$ such that
\begin{equation}
\intO \nabla u\cdot\nabla (v-u) \textrm{d}x\geq\intO f(v-u)\textrm{d}x \quad \mbox{for all }v \in K. \label{2varineq}
\end{equation}
The convex constrain $v\in K$ can be transformed into a linear term containing a new (Lagrange) variable in

\begin{problem}[Perturbed problem] \label{prob2} Let $W:=\{v+t\phi:v\in V_{0}\; \textrm{and}\; t\in\mathbb{R}\}\subset V$.
For given
\begin{equation} \label{nezapor}
\mu \in \Lambda:=\{\mu \in W^{*}: \left< \mu,v-\phi\right>\geq0 \;\;
\mbox{for all }\; v\in K\}
\end{equation}
find $u_{\mu}\in V_{0}$ such that
\begin{equation} \label{2defpert}
J_{\mu}(u_{\mu})=\inf\limits_{v\in V_{0}} J_{\mu}(v),
\end{equation}
where
$\left< \cdot,\cdot \right>$ denote the duality pairing of $W$ and $W^{*}$ and the perturbed functional $J_{\mu}$ is defined as
\begin{equation}  \label{2Pert}
J_{\mu}(v):=J(v)-\left< \mu,v-\phi \right>.
\end{equation}
\end{problem}

\noindent Problems \ref{prob1} and \ref{prob2} are related and it obviously holds
\begin{equation} \label{2ineqJmuJ}
J_{\mu}(u_\mu) \leq J(u) \quad \mbox{for all }\mu \in \Lambda.
\end{equation}

\begin{lemma}[Existence of optimal multiplier]
There exists
$\lambda \in \Lambda$ such that
\begin{equation}  \label{prvek}
u_{\lambda}=u
\end{equation}
and
\begin{equation} \label{hodnota}
J_{\lambda}(u)= J(u).
\end{equation}
\begin{proof}
Let $w\in W$ is arbitrary. We decompose
\begin{equation} \label{rozklad}
w=v+t\phi,
\end{equation}
where $v\in V_{0}$ and $t\in\mathbb{R}$ and this decomposition can be shown to be unique. Now,
we define a functional $\lambda$ as follows:
\begin{equation} \label{fc}
\left\langle\lambda,w\right\rangle :=\intO \nabla u\cdot\nabla v \,\textrm{d}x -\intO  f v \textrm{d}x+t\left[\intO \nabla u\cdot\nabla u \,\textrm{d}x -\intO  f u \textrm{d}x\right].
\end{equation}
We assert that the functional defined by (\ref{fc}) has required properties (\ref{prvek}) and (\ref{hodnota}). Apparently, $\lambda$ is a linear functional on $W$.
The functional $\lambda$ is also continuous. It is a consequence of continuity of decomposition (\ref{rozklad}), which can be proved as follows.
Let $w \in W$ is arbitrary and $w_{n}\rightarrow w$ in $W$, where $w_{n}\in W$. With respect of (\ref{rozklad}), we can write $w_{n}=v_{n}+t_{n}\phi$ and $w=v+t\phi$, where $v_{n}, v\in V_{0}$ and $t_{n},t\in\mathbb{R}$. If we use the unique orthogonal decomposition of element  $\phi\in V$, we infer that
\begin{equation}
|t_{n}-t|\|\phi^{\bot}\|_{V}\leq \|w_{n}-w\|_{V},   \label{n1}
\end{equation}
where
$\phi^{\bot}$ is the component of $\phi$ orthogonal to subspace $V_{0}$. Moreover, it follows from the triangle inequality that
\begin{equation}
\|v_{n}-v\|_{V}\leq \|w_{n}-w\|_{V}+|t_{n}-t|\|\phi\|_{V}.  \label{n2}
\end{equation}
As a consequence of (\ref{n1}) and (\ref{n2}), $t_{n}\rightarrow t$ and $v_{n}\rightarrow v$ in $W$. Thus, the decomposition (\ref{rozklad})  is continuous.
Now, if we restrict the space $W$ to the origin $V_{0}$, we obtain
\begin{equation} \label{fcV}
\left\langle\lambda,w\right\rangle :=\intO \nabla u\cdot\nabla w \,\textrm{d}x -\intO  f w \textrm{d}x\quad \mbox{for all } w\in V_0,
\end{equation}
which is equivalent to $\inf\limits_{w\in V_{0}} J_{\lambda}(w)=J_{\lambda}(u)$, i.e.,  the property (\ref{prvek}) is fulfilled.
Furthermore, if we take $w=u-\phi$, it follows from (\ref{fc}) that
$$
\left\langle\lambda,u-\phi\right\rangle=\left\langle\lambda,u\right\rangle-\left\langle\lambda,\phi\right\rangle=0
$$
and consequently the property (\ref{hodnota}) is fulfilled. Finally, we should verify the condition of nonnegativity from definition (\ref{nezapor}). Let $v\in K$
is arbitrary. It follows from (\ref{2varineq}) and (\ref{fc}) that
$$
\left< \lambda,v-\phi\right>=\left< \lambda,v-u\right>\geq0.
$$
\end{proof}

\end{lemma}

\begin{remark}[Existence of optimal multiplier in the case of nonpositive obstacle $\phi\in V_{0}$]
If we would deal with a nonpositive obstacle $\phi\in V_{0}$, the existence of optimal multiplier could be proved as follows.
Once again, the relation (\ref{fcV}) defines a linear continuous functional $\lambda$ in $V_{0}$ such that $u$ minimizes the perturbed functional $J_{\mu}$ defined by (\ref{2Pert}) with $\mu =\lambda$.
Since $\phi\in K$, we can  apply the inequality (\ref{2varineq}) to $v=\phi$ and $v=2u-\phi$. Consequently, we obtain that
$
\left\langle\lambda,u-\phi\right\rangle
=0.$
Subsequently, it follows from (\ref{2Pert}) that
the property (\ref{hodnota}) is fulfilled. The condition of nonnegativity from definition (\ref{nezapor}) is also fulfilled. It follows from the inequality (\ref{2varineq}) if we put
$v=u+w$, where $w\in V_{0}$, $w\geq0$ a.e. in $(\Omega)$.
\end{remark}

\begin{remark}[Representation of (\ref{fc}) by a nonnegative function $\lambda\in L^{2}(\Omega)$]
If $u$ has a higher regularity,
\begin{equation} \label{u_H2}
u\in V_{0}\cap H^{2}(\Omega),
\end{equation} then integration by parts yields
\begin{equation} \label{lambda_L2}
\left\langle\lambda,w\right\rangle =-\intO \Delta u \,v\textrm{d}x -\intO  f v \textrm{d}x+t\left[-\intO \Delta u \,u\textrm{d}x -\intO  f u \textrm{d}x\right]
=\intO \lambda v\textrm{d}x+t\intO \lambda u\textrm{d}x
\end{equation}
for all $w \in W$, where
\begin{equation} \label{podmlambda}
\lambda =-(\Delta u \,+f).
\end{equation}
We show additionaly that
\begin{equation}
\label{lambda_nonnegative}
\lambda\geq0 \mbox{ a.e. in } \Omega
\end{equation} by choosing $w \in V_0, w\geq0$ a.e. in $\Omega$. Then $v: = u + w \in K$ and inequality  (\ref{2varineq}) rewrites as
$$
\intO \lambda w\textrm{d}x=\intO \nabla u\cdot\nabla w \,\textrm{d}x -\intO  f w \textrm{d}x\geq0,
$$
which implies \eqref{lambda_nonnegative}.
\label{remark_lambda_in_L2}
\end{remark}

\section{Functional a posteriori error estimate}
\noindent We are interested in analysis and numerical properties of the a posteriori error estimate in the energetic norm
$$
\|v\|_{E}:=\left(\intO\nabla v\cdot\nabla v\,\textrm{d}x\right)^{\frac{1}{2}}.
$$
This section is based on results of S. Repin et al. \cite{BuRe, NeRep, Re4}.
It is simple to see that
\begin{eqnarray}\label{energy_difference}
J(v)-J(u)=\frac{1}{2}\intO\nabla (v-u)\cdot\nabla (v-u)\,\textrm{d}x+\intO \nabla u\cdot\nabla (v-u)\,\textrm{d}x-\intO f(v-u)\textrm{d}x 
\end{eqnarray}
for all $v\in K$
and \eqref{2varineq} implies the energy estimate
\begin{equation}\label{energyEstimate}
\frac{1}{2}\|v-u\|^{2}_{E}\leq J(v)-J(u) \quad \mbox{for all } v\in K.
\end{equation}

\begin{remark}[Sharpness of estimate \eqref{energyEstimate}] \label{remark:sharpness_2D}
It is clear from \eqref{energy_difference}, the estimate \eqref{energyEstimate} turns into equality if
\begin{equation}\label{estimate_equality_condition_2D}
\left< \lambda,  v-u \right>  = \intO \nabla u\cdot\nabla (v-u)\,\textrm{d}x- \intO f(v-u)\textrm{d}x =0 \quad \mbox{for all } v\in K.
\end{equation}
This situation always occurs if $\lambda=0$. Then, \eqref{fcV} implies that $u$ is a solution of Problem \ref{prob1} in the whole space $V_0$. This corresponds to a linear problem without any obstacle. However, the estimate \eqref{energyEstimate} can turn into equality also for the active obstacle.
We discuss it further in Section \ref{sec:exact_solution}.
\end{remark}

Estimate \eqref{energyEstimate} can only be tested for problems with known exact solution $u\in K$.
By using \eqref{2ineqJmuJ}, we obtain the estimate
\begin{equation}\label{eq9}
J(v)-J(u) \leq J(v)-J_{\mu}(u_{\mu})  \quad \mbox{for all }\mu \in \Lambda.
\end{equation}
In practical computations, $u_{\mu} \in V_{0}$ will be approximated by $u_{\mu,h} \in V_{0,h}$ from some finite dimensional subspace $ V_{0,h} \subset V_{0}$ (see Section \ref{sec:numerics} for details). Therefore, it holds
$$ J_{\mu}(u_{\mu,h}) \geq J_{\mu}(u_\mu) $$
and $J_{\mu}(u_{\mu,h})$ can not replace $J(u)$ in \eqref{eq9} so that the inequality holds.
To avoid this difficulty, we establish the following dual problem:

%
\begin{problem}[Dual perturbed problem] Find $\tau^{*}_{\mu}\in Q^{*}_{f\mu}\subset [L^{2}(\Omega)]^2$ such that
\begin{equation}
J^{*}_{\mu}(\tau^{*}_{\mu})=\sup\limits_{q^{*}\in Q^{*}_{f\mu}} J^{*}_{\mu}(q^{*}),
\end{equation}
where
\begin{equation}  \label{dupert}
J^{*}_{\mu}(q^{*})=-\frac{1}{2}\intO q^{*}\cdot q^{*} \textrm{d}x + \left< \mu, \phi \right>
\end{equation}
and
\begin{equation}   \label{mnozQ}
Q^{*}_{f\mu}:=\left\{q^{*}\in [L^{2}(\Omega)]^2 :\left< \mu, v \right> = \intO q^{*}\cdot \nabla v\,\textrm{d}x - \intO f v\textrm{d}x
\quad \mbox{for all } v\in V_{0}\right\}.
\end{equation}
\end{problem}

\begin{lemma}
It holds
$$
\sup\limits_{q^{*}\in Q^{*}_{f\mu}} J^{*}_{\mu}(q^{*})=J^{*}_{\mu}(\nabla u_{\mu})=J_{\mu}(u_\mu).
$$
\begin{proof}
As a consequence of (\ref{2defpert}), it holds that $\nabla u_{\mu}\in Q^{*}_{f\mu}$. Let $w\in Q^{*}_{f\mu}$
is arbitrary. Since
$$
J^{*}_{\mu}(w)=J^{*}_{\mu}(\nabla u_{\mu})-\intO \nabla u_{\mu}\cdot(w-\nabla u_{\mu})\textrm{d}x-\frac{1}{2}\intO (w-\nabla u_{\mu})\cdot(w-\nabla u_{\mu})\textrm{d}x
$$
and\;
$
\intO \nabla u_{\mu}\cdot(w-\nabla u_{\mu})\textrm{d}x=0
$
\;in consequence of \;$\nabla u_{\mu},w\in Q^{*}_{f\mu}$,
we deduce that $J^{*}_{\mu}(\nabla u_{\mu})$ is supremum of dual perturbed functional $J^{*}_{\mu}$. Finally, it is not difficult to verify that
$J^{*}_{\mu}(\nabla u_{\mu})=J_{\mu}(u_\mu)$.

\end{proof} \label{lemma2}
\end{lemma}
\begin{corollary}
If we choose
$$
\mu=\lambda,
$$
it holds
$$
J_{\lambda}(u)=\inf\limits_{v\in V_{0}} J_{\lambda}(v)=\sup\limits_{q^{*}\in Q^{*}_{f\lambda}} J^{*}_{\lambda}(q^{*})=J^{*}_{\lambda}(\nabla u)=J(u).
$$
\end{corollary}
\bigskip
\noindent It follows from Lemma \ref{lemma2}, we can replace inequality (\ref{eq9}) by
\begin{equation}
J(v)-J(u)\leq J(v)-\sup\limits_{q^{*}\in Q^{*}_{f\mu}} J^{*}_{\mu}(q^{*})
\leq J(v)-J^{*}_{\mu}(q^{*}),      \label{eq11}
\end{equation}
where $q^{*}\in Q^{*}_{f\mu}$ is arbitrary.
The practical limitation of estimate \eqref{eq11} is to satisfy the constrain
$q^{*}\in Q^{*}_{f\mu}$.
From now, we consider a special case of the multiplier defined as
\begin{equation}\label{defFun}
\left< \mu, w \right> := \intO \mu w \,\textrm{d}x,
\end{equation}
where
\begin{equation}\label{admis}
\mu \in \Lambda:=\left\{\mu\in L^2(\Omega): \mu\geq0\; \mbox{a.e.}\; \textrm{in}\; \Omega\right\}.
\end{equation}
S. Repin transformed \eqref{eq11} in the so called majorant estimate

\begin{equation}\label{majorantEstimate}
J(v)-J(u)\leq \mathcal{M}(v,f,\phi;\beta,\mu,\tau^{*}),
\end{equation}
where the right-hand side of \eqref{majorantEstimate} denotes the functional majorant
\begin{multline} \label{majorantForm}
\mathcal{M}(v,f,\phi;\beta,\mu, \tau^{*}) :=
\frac{1+\beta}{2}\intO (\nabla v-\tau^{*})\cdot(\nabla v-\tau^{*}) \textrm{d}x \\
+\frac{1}{2}\left(1+\frac{1}{\beta}\right) C_{\Omega}^2 \|{\rm{div}} \,\tau^{*}+f+\mu\|^{2}_{L^{2}(\Omega)}
+\intO \mu(v-\phi)\textrm{d}x,
\end{multline}
where a constant $C_{\Omega}>0$ originates from the Friedrichs inequality
$$
\intO u^{2}\textrm{d}x\leq C_{\Omega}^{2} \intO \nabla u\cdot\nabla u \textrm{d}x\quad\forall u\in V_{0}.
$$
Estimate \eqref{majorantEstimate} is valid for $\beta>0$, $\mu \in \Lambda$ and $\tau^{*}\in H(\Omega,{\rm{div}})$,
where $$H(\Omega,{\rm{div}}):=\{\tau^{*}\in [L^2(\Omega)]^2:{\rm{div}}\,\tau^{*}\in L^{2}(\Omega)\}.$$

\begin{lemma}[Optimal majorant parameters] \label{rem:optimal}
Suppose
(\ref{defFun}) - (\ref{admis}) and, let the assumption (\ref{u_H2}) is fulfilled.
If we choose $\tau^{*}=\nabla u$, $\mu=\lambda\in L^{2}(\Omega)$ and $\beta\rightarrow 0$, then,
the inequality in (\ref{majorantEstimate}) changes to equality, i.e.
the majorant on right-hand side of (\ref{majorantEstimate}) defines the difference of energies\, $J(v)-J(u)$\, exactly.
\begin{proof}
If $\mu=\lambda$ and $\tau^{*}=\nabla u$, it is consequence of (\ref{podmlambda}) that the second term on the right-hand side of  (\ref{majorantForm}) vanishes.
Moreover, the last term can be written as
$$
\intO \lambda(v-\phi)\textrm{d}x=\intO \lambda(v-u)\textrm{d}x=\intO \nabla u\cdot\nabla(v-u)\textrm{d}x-\intO f(v-u)\textrm{d}x
$$
and consequently, if $\beta=0$, it follows from (\ref{energy_difference}) that the majorant with optimal parameters estimates the difference of energies $J(v)$ and $J(u)$\, exactly.
\end{proof}
\end{lemma}

\noindent Practically, optimal parameters are unknown. For given solution approximation $v$, loading $f$ and the obstacle $\phi$, the majorant $\mathcal{M}$ represents a convex functional in each of variables $\beta$, $\mu$ ,$\tau^{*}$.  Our goal is to find, at least approximately,
such variables $\beta_{\mathrm{opt}}$, $\mu_{\mathrm{opt}}$ and $\tau^{*}_{\mathrm{opt}}$ that minimize the majorant $\mathcal{M}$.

\begin{problem}[Majorant minimization problem] Let $v\in K$, $f\in L^2(\Omega)$, $\phi<0$ be given. Find optimal $\beta_{\mathrm{opt}}>0$, $\mu_{\mathrm{opt}}\in\Lambda$ and $\tau^{*}_{\mathrm{opt}}\in H(\Omega,{\rm{div}})$ such that
$$
(\beta_{\mathrm{opt}},\mu_{\mathrm{opt}},\tau^{*}_{\mathrm{opt}}) = \argmin \limits_{\beta,\mu,\tau^{*}}\mathcal{M}(v,f,\phi;\beta,\mu,\tau^{*}).
$$
\end{problem}

\noindent To this end, we use the following minimization algorithm:

\begin{algorithm}[Majorant minimization algorithm]\label{Alg:Majorant}
\noindent Let $k=0$ and let $\beta_{k}>0$ and $\mu_{k}\in\Lambda$ be given. Then:
\begin{itemize}
\item[(i)] find $\tau^{*}_{k+1}\in H(\Omega,{\rm{div}})$ such that
$$
\tau^{*}_{k+1}=\argmin\limits_{\tau^{*}\in H(\Omega,{\rm{div}})}\mathcal{M}(v,f,\phi;\beta_k,\mu_k,\tau^{*}),
$$
\item[(ii)] find $\mu_{k+1}\in \Lambda$ such that
$$
\mu_{k+1}=\argmin\limits_{\mu\in \Lambda}\mathcal{M}(v,f,\phi;\beta_k,\mu,\tau^{*}_{k+1}),
$$
\item[(iii)] find $\beta_{k+1}>0$ such that
$$
\beta_{k+1}=\argmin\limits_{\beta > 0 }\mathcal{M}(v,f,\phi;\beta,\mu_{k+1},\tau^{*}_{k+1}),
$$
\item[(iv)] set $k=k+1$ are repeat (i)-(iii) until convergence.
\end{itemize}
\end{algorithm}

\begin{remark}[Functional majorant in 1D] \label{rem3}
The goal is to verify the majorant error estimate for obstacle problem in 1D. In this simplified case, the former domain $\Omega$ reduces
to one-dimensional interval $(0,1)$. We set $V:=H^{1}(0,1)$ and $V_{0}:=H^{1}_{0}(0,1)$, the energy functional (\ref{2rov2}) reads
$$
J(v) = \frac{1}{2}\int_0^1(v')^{2}\textrm{d}x-\int_0^1fv\,\textrm{d}x
$$
and the admissible set is defined as
$$
K:=\bigl\{v \in V_{0}:\,
v(x) \geq \phi(x) \mbox{ a.e.}\;\textrm{in}\; (0,1) \},
$$
where $f\in L^{2}(0,1)$ and $\phi\in V$ such that $\phi\not\in V_{0}$ and $\phi(x)<0 \mbox{ a.e.}\;\textrm{in}\; (0,1)$.
Then, the functional majorant $\mathcal{M}$ takes the form
\begin{equation}
\mathcal{M}(v,f,\phi;\beta,\mu,\tau^{*})=\frac{1+\beta}{2}\int_{0}^{1}(v'-\tau^{*})^{2}\textrm{d}x
+\frac{1}{2}\left(1+\frac{1}{\beta}\right)\|(\tau^{*})'+f+\mu\|^{2}_{L^{2}(0,1)}
+\int_{0}^{1}\mu(v-\phi)\textrm{d}x,
\end{equation}
where $\beta>0$, $\tau^{*}\in V$ and $\mu \in \Lambda =\left\{\mu\in L^2(0,1): \mu\geq0\; \mbox{a.e.}\; \textrm{in}\; (0,1)\right\}$.
\end{remark}

\begin{remark}[Majorant minimization in 1D]
In 1D case, the minimization in step (i) is equivalent to the
following variational equation
: Find $\tau^{*}_{k+1} \in V$ such that
\begin{equation} \label{minization1D_step_i}
\begin{aligned}
(1+\beta_{k})\int_{0}^{1}\tau^{*}_{k+1} w\textrm{d}x
&+\left(1+\frac{1}{\beta_{k}}\right)\int_{0}^{1}(\tau^{*}_{k+1})'w'\textrm{d}x\\
&=(1+\beta_{k})\int_{0}^{1}v'w\textrm{d}x
-\left(1+\frac{1}{\beta_{k}}\right)\int_{0}^{1}(f+\mu_k)w'\textrm{d}x
\qquad \mbox{for all } w\in V.
\end{aligned}
\end{equation}
The minimization in step (ii) is equivalent to the variational inequality:
Find $\mu_{k+1}\in \Lambda$ such that
\begin{equation}  \label{minization1D_step_ii}
\int_{0}^{1}\left[\left(1+\frac{1}{\beta_{k}}\right)\left[\mu_{k+1}+(\tau^{*}_{k+1})'+f\right]+v-\phi\right](w-\mu_{k+1})\textrm{d}x\geq0
\quad \mbox{for all } w\in \Lambda.
\end{equation}
The minimization in step (iii) leads to the explicit relation
\begin{equation}  \label{minization1D_step_iii}
\beta_{k+1}=\frac{\|(\tau^{*}_{k+1})'+f+\mu_{k+1}\|_{L^{2}(0,1)}}{\|v'-\tau^{*}_{k+1}\|_{L^{2}(0,1)}}.
\end{equation}
\end{remark}

\section{1D benchmark with known analytical solution}\label{sec:exact_solution}
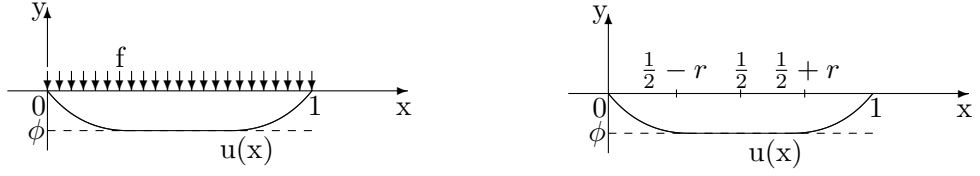
\begin{figure}
\begin{minipage}{9cm}
\setlength{\unitlength}{1.5pt}
\begin{picture}(45,60)(-60,0)
\put(-10,30){\vector(1,0){100}}
\put(0,36){\vector(0,1){15}}
\put(0,15){\line(0,1){15}}
\put(0,35){\vector(0,-1){5}}
\put(3,35){\vector(0,-1){5}}
\put(6,35){\vector(0,-1){5}}
\put(9,35){\vector(0,-1){5}}
\put(12,35){\vector(0,-1){5}}
\put(15,35){\vector(0,-1){5}}
\put(18,35){\vector(0,-1){5}}
\put(21,35){\vector(0,-1){5}}
\put(24,35){\vector(0,-1){5}}
\put(27,35){\vector(0,-1){5}}
\put(30,35){\vector(0,-1){5}}
\put(33,35){\vector(0,-1){5}}
\put(36,35){\vector(0,-1){5}}
\put(39,35){\vector(0,-1){5}}
\put(42,35){\vector(0,-1){5}}
\put(45,35){\vector(0,-1){5}}
\put(48,35){\vector(0,-1){5}}
\put(51,35){\vector(0,-1){5}}
\put(54,35){\vector(0,-1){5}}
\put(57,35){\vector(0,-1){5}}
\put(60,35){\vector(0,-1){5}}
\put(63,35){\vector(0,-1){5}}
\put(66,35){\vector(0,-1){5}}

\multiput(0,20)(4,0){17}{\line(1,0){2}}

\put(87,24){x}
\put(65,24){1}
\put(17,36){f}
\put(-4,49){y}
\put(43,13){u(x)}
\put(-5,18){$\phi$}
\qbezier(0,30)(8,20)(20,20)
\qbezier(46,20)(58,20)(66,30)
\put(20,20){\line(1,0){26}}
\put(-4,24){0}
\end{picture}
\end{minipage}

\begin{minipage}{9cm}
\setlength{\unitlength}{1.5pt}
\begin{picture}(0,0)(-200,0)
\put(-10,30){\vector(1,0){100}}
\put(0,15){\vector(0,1){35}}

\multiput(0,20)(4,0){17}{\line(1,0){2}}

\put(87,24){x}
\put(65,24){1}
\put(-4,49){y}
\put(35,13){u(x)}
\put(-5,18){$\phi$}
\qbezier(0,30)(8,20)(20,20)
\qbezier(46,20)(58,20)(66,30)
\put(20,20){\line(1,0){26}}
\put(-4,24){0}
\put(31,34){$\frac{1}{2}$}
\put(8,34){$\frac{1}{2}-r$}
\put(41,34){$\frac{1}{2}+r$}

\put(17,29){\line(0,1){2}}

\put(33,29){\line(0,1){2}}

\put(49,29){\line(0,1){2}}
\end{picture}
\end{minipage}
\caption{Benchmark setup: Constant forces $f$ pressing continuum against a constant lower obstacle $\phi$, exact displacement $u$ (left) and construction of exact displacement $u$ in detail (right).}
\label{fig:benchmark}
\end{figure}

\begin{figure}
\center
\includegraphics[width=0.7\textwidth]{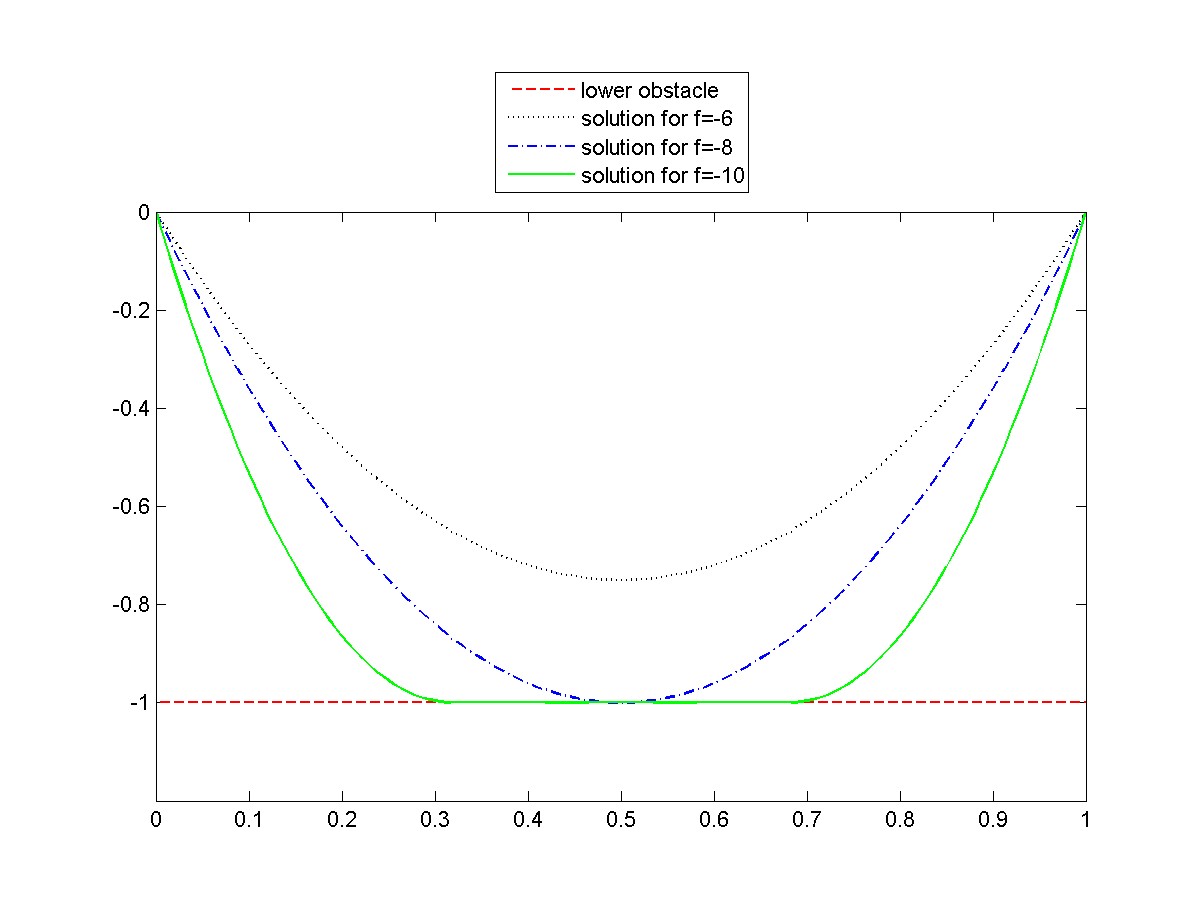}
\caption{Solutions for problems with loadings $f \in \{-6, -8, -10 \}$ and $\phi=-1$.}\label{Fig:solution_f}
\end{figure}

We derive an exact solution of Problem \ref{prob1} - modified to 1D problem (see Remark \ref{rem3}) - assuming negative constant functions $f$ and $\phi$ . The resulting solution is displayed in Figure \ref{fig:benchmark}
for the case of active obstacle. A mechanical intuition suggests that for small values (considered in absolute value) of acting force $f$, there will be no contact with the obstacle and there will be a contact on a subset of interval $(0, 1)$ located symmetrically around the value $x=1/2$ for higher values of $f$.
The solution of Problem \ref{prob1} with inactive obstacle reads
\begin{equation}
u(x) = \frac{f}{2}(x - x^2).        \label{eq4}
\end{equation}
The minimal value of $u$ on interval (0,1) is attained at $x=1/2$ and the inactive obstacle condition $u(1/2) > \phi$ is satisfied for
\begin{equation}\label{condition_inactive_obstacle}
|f| < 8|\phi|.
\end{equation}
Then, the corresponding energy reads
\begin{equation}
J(u) := -\frac{f^2}{24}.    \label{eq5}
\end{equation}
\noindent The obstacle is active if
\begin{equation}\label{condition_active_obstacle}
|f| \geq 8|\phi|,
\end{equation}
and the solution has the following form
\begin{equation*}
u(x)=\left\{
\begin{array}{lrl}
\medskip -\frac{f}{2}x^{2}+\frac{\phi+\frac{f}{2}(\frac{1}{2}-r)^{2}}{\frac{1}{2}-r}x & \quad\textrm{if} & x\in [0,\frac{1}{2}-r) \\
\medskip \phi & \quad\textrm{if} & x \in [\frac{1}{2}-r,\frac{1}{2}+r] \\
\medskip -\frac{f}{2}(x-1)^{2}-\frac{\phi+\frac{f}{2}(\frac{1}{2}-r)^{2}}{\frac{1}{2}-r}(x-1) & \quad\textrm{if} & x\in (\frac{1}{2}+r,1]
\end{array} \right.
\end{equation*}
for unknown parameter $r\in [0,\frac{1}{2}]$. The parameter $r$ determines the active contact set $[\frac{1}{2}-r, \frac{1}{2}+r]$ and its value can be determined from the
minimum of energy
\begin{equation*}
J(u) =\frac{[\phi+\frac{f}{2}(\frac{1}{2}-r)^{2}]^2}{\frac{1}{2}-r}
-2[\phi+\frac{f}{2}(\frac{1}{2}-r)^{2}]f(\frac{1}{2}-r)    \label{eq7}
+\frac{2f^2}{3}(\frac{1}{2}-r)^{3}-2fr\phi
\end{equation*}
over all value of $r\in [0,\frac{1}{2}]$. The minimal energy
\begin{equation}
J(u)=f \phi (\frac{4}{3} \sqrt{\frac{2\phi}{f}}-1)
\end{equation}
is achieved for the argument
\begin{equation}
r=\frac{1}{2}-\sqrt{\frac{2\phi}{f}}.     \label{eq8}
\end{equation}
Therefore, the solution of the problem with the active obstacle reads
\begin{equation} \label{soll}
u(x)=\left\{
\begin{array}{lrl}
\medskip -\frac{f}{2}x^{2}-\sqrt{2\phi f}x & \quad\textrm{if} & x\in \left[0,\sqrt{\frac{2\phi}{f}}\right) \\
\medskip \phi & \quad\textrm{if} & x \in \left[\sqrt{\frac{2\phi}{f}},1-\sqrt{\frac{2\phi}{f}}\right] \\
\medskip -\frac{f}{2}(x-1)^{2}+\sqrt{2\phi f}(x-1) & \quad\textrm{if} & x\in \left(1-\sqrt{\frac{2\phi}{f}},1\right]
\end{array} \right.
\end{equation}
Figure \ref{Fig:solution_f} provides few numerical approximations of $u$, see Section \ref{sec:numerics} for details.
The first-order derivative
\begin{equation}
u'(x)=\left\{
\begin{array}{lrl}
\medskip -fx-\sqrt{2\phi f} & \quad\textrm{if} & x\in \left[0,\sqrt{\frac{2\phi}{f}}\right) \\
\medskip 0 & \quad\textrm{if} & x \in \left[\sqrt{\frac{2\phi}{f}},1-\sqrt{\frac{2\phi}{f}}\right] \\
\medskip -f(x-1)+\sqrt{2\phi f} & \quad\textrm{if} & x\in \left(1-\sqrt{\frac{2\phi}{f}},1\right]
\end{array} \right.
\end{equation}
is continuous everywhere. It is not difficult to show that\begin{equation}
u''(x)=\left\{
\begin{array}{lrl}
\medskip -f & \quad\textrm{if} & x\in \left(0,\sqrt{\frac{2\phi}{f}}\right) \\
\medskip 0 & \quad\textrm{if} & x \in \left(\sqrt{\frac{2\phi}{f}},1-\sqrt{\frac{2\phi}{f}}\right) \\
\medskip -f & \quad\textrm{if} & x\in \left(1-\sqrt{\frac{2\phi}{f}},1\right)
\end{array} \right.
\end{equation}
is the second-order weak derivative of $(\ref{soll})$. With respect to (\ref{podmlambda}),
the optimal multiplier
for our 1D benchmark problem reads
\begin{equation}\label{lambda_exact}
\lambda(x)=\left\{
\begin{array}{lrl}
\medskip 0 & \quad\textrm{if} & x\in \left(0,\sqrt{\frac{2\phi}{f}}\right) \\
\medskip -f & \quad\textrm{if} & x \in \left(\sqrt{\frac{2\phi}{f}},1-\sqrt{\frac{2\phi}{f}}\right) \\
\medskip 0 & \quad\textrm{if} & x\in \left(1-\sqrt{\frac{2\phi}{f}},1\right)
\end{array} \right.
\end{equation}
so that it is a piecewise constant function.

\begin{remark}[Sharpness of estimate \eqref{energyEstimate} for 1D benchmark] \label{remark:sharpness}
It is easy to show that the estimate \eqref{energyEstimate} can turn into equality for the active obstacle.
Indeed, in our 1D benchmark, the condition \eqref{estimate_equality_condition_2D} rewrites as
\begin{equation}
\left< \lambda,  v-u \right>  =
-f\int_{\sqrt{\frac{2\phi}{f}}}^{1-\sqrt{\frac{2\phi}{f}}} (v-u)\textrm{d}x = 0,  \label{ner2}
\end{equation}
if the contact zone of an approximate solution $v\in K$ includes
whole contact zone $\left[\sqrt{\frac{2\phi}{f}},1-\sqrt{\frac{2\phi}{f}}\right]$ of exact solution $u$.
\end{remark}

\section{Numerical experiments}\label{sec:numerics}
A MATLAB software is available as a package {\it Obstacle problem in 1D and its a posteriori error estimate} at Matlab Central under
{\small\url{http://www.mathworks.com/matlabcentral/fileexchange/authors/37756}}. \\

\begin{figure}
\center
\includegraphics[width=0.8\textwidth]{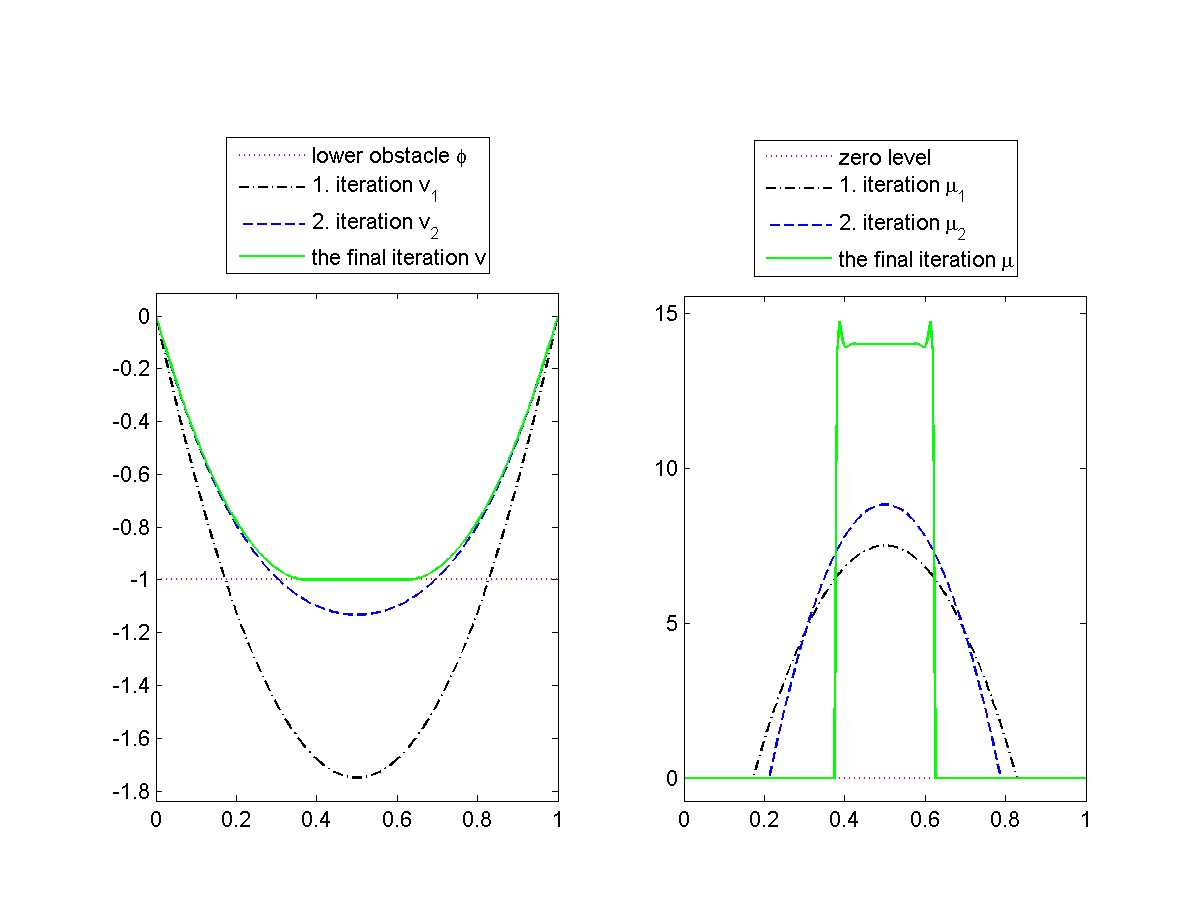}
\caption{The first, the second and the final (the 10000th) iterations of the Uzawa algorithm run on an uniform mesh with $641$ nodes for the loading $f=-14$ and the obstacle $\phi=-1$.}\label{Fig:Uzawa}
\includegraphics[width=0.8\textwidth]{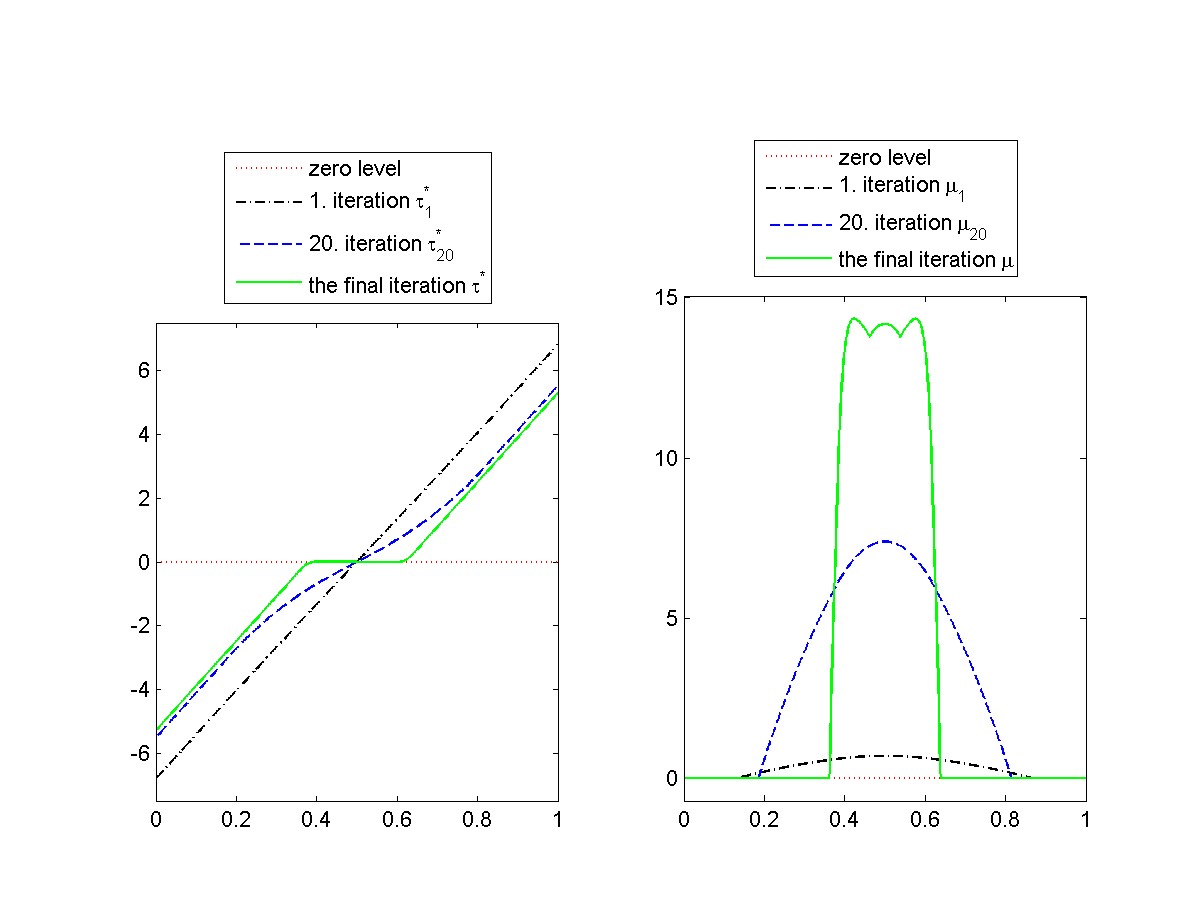}
\caption{The first, the twentieth and the final (the 10000th) iterations of the majorant minimization algorithm run on an uniform mesh with $641$ nodes for the loading $f=-14$ and the obstacle $\phi=-1$. We assumed the initial setup $\beta_0=1, \mu_0=0$ and the approximation $v$ obtained after 100 iterations of the Uzawa algorithm.}\label{Fig:Majorant}
\end{figure}

Assuming the interval partition $\mathcal{T}$ with $n$ nodes
\begin{equation} \label{interval_partition}
0=x_{1}<x_{2}< \ldots <x_{n}=1,
\end{equation}
we define $V_h \subset V$ as the finite dimensional space of nodal linear functions with a basis $\psi_{j}$, $j=1\hdots n$ and its subspace $V_{0,h}$ of functions satisfying homogeneous Dirichlet boundary conditions.
Using these basis functions, a stiffness matrix $A=(a_{ij})$ and a mass matrix $M=(m_{ij})$ are defined as
$$
a_{ij}:=\int_{0}^{1}\psi'_{i}\psi'_{j}\textrm{d}x  \quad m_{ij}=\int_{0}^{1}\psi_{i}\psi_{j}\textrm{d}x.
$$
A numerical approximation $v \in V_{0,h}$ of the exact solution $u \in V_0$ is constructed by the Uzawa algorithm.
\begin{algorithm}[Uzawa algorithm] \label{Alg:Uzawa}
\noindent
\begin{enumerate}
	\item Set the initial Lagrange multiplier $\mu_0=0$.
	\item Start of the loop: for $k=1,2, \dots$ do until convergence:
	\item Find an approximation $v_{k} \in  V_{0,h}$ such that $J_{\mu_k}(v_k) \rightarrow \min$.
	\item Set a Lagrange multiplier $\mu_k=(\mu_{k-1} + \rho (v_k - \phi))^+$.
	\item End of the loop.
	\item Output $v=v_k$ and $\mu =  \mu_k$.
\end{enumerate}
\end{algorithm}
\noindent
The approximation $v_{k} =\sum_{j=2}^{n-1}v_{k,j}\psi_{j}$ in step 3. of Algorithm \ref{Alg:Uzawa} is computed from
the equivalent variational equation
$$
\int_{0}^{1}v'_{k} w'\textrm{d}x=\int_{0}^{1}(f+\mu_{k}) w \, \textrm{d}x\qquad \mbox{for all } w\in V_{0,h}
$$
leading to a linear system of equations for coefficients $v_{k,2},\dots,v_{k,n-1}$. The convergence of Algorithm \ref{Alg:Uzawa} depends on the choice of the scalar parameter $\rho$ and it can be shown, see e.g. \cite{GLT}, that is alway converges for $\rho \in (0, \rho_1)$ for some $\rho_1 > 0$. Some iterations of Algorithm \ref{Alg:Uzawa} with $\rho=10$ are displayed in Figure \ref{Fig:Uzawa}. Algorithm \ref{Alg:Uzawa} converges slowly and therefore lower number of its iterations provides a poor approximation $v$ of the exact solution $u$. In the following, we consider three particular sets of approximations $v$ obtained by Algorithm \ref{Alg:Uzawa} with different numbers of iterations:
\begin{description}
	\item \quad a) 100 iterations,
\quad b) 1000 iterations,
\quad c) 10000 iterations.
\end{description}
The sets of solutions a), b), c) will be constructed for the uniform mesh $\mathcal{T}$ with $641$ nodes (which corresponds to 6 uniform refinements of an initial uniform mesh with 10 elements) and for various loadings
$$f \in \{-5, -6, \dots, , -17, -18\}.$$
It follows from \eqref{condition_inactive_obstacle} and \eqref{condition_active_obstacle},  the obstacle is
inactive for $f \in \{-5, \dots, -7\}$ and active for $f \in \{-8, \dots, -18\}$. Therefore, Algorithm \ref{Alg:Uzawa} converges in a continuous setup for $f \in \{-5, \dots, -7\}$ after one iteration and approximations a), b), c) coincide. A verification of the energy estimate \eqref{energyEstimate} is reported in Tables \ref{tab:diferent_f_HundredIterations}, \ref{tab:diferent_f_ThousandIterations}, \ref{tab:diferent_f_tenThousandIterations}. We notice that the gap between the energy error
$\frac{1}{2}\|v-u\|^{2}_{E}$ and the difference of energies $J(v)-J(u)$ is very small for approximations c) and becomes larger for approximations b) and a). In the case of inactive contact, the gap is apparently zero, see Remark \ref{remark:sharpness_2D}.

For the verification of the majorant estimate \eqref{majorantEstimate}, we run a discretized version of Algorithm \ref{Alg:Majorant}.
The minimal argument $\tau^{*}_{k+1} \in V_h$ in step (i) of Algorithm \ref{Alg:Majorant} is searched in the form $\tau^{*}_{k+1}=\sum_{j=1}^{n}y_{j}\psi_{j}$, where
coefficients $y=(y_1,\dots,y_n) \in R^n$ follow (see \eqref{minization1D_step_i}) from a linear system of equations
\begin{equation}
\left[(1+\beta_k)M+\left(1+\frac{1}{\beta_k}\right)A\right]y=(1+\beta_k)b-\left(1+\frac{1}{\beta_k}\right)c,
\label{eq:system_y}
\end{equation}
where $b$ and $c$ are $n$-\;dimensional vectors defined as
$$
b_i=\int_{0}^{1}v'\psi_{i}\textrm{d}x, \quad c_i=\int_{0}^{1}(f+\mu_k)\psi'_{i}\textrm{d}x
$$
for $i,j=1\hdots n$. The minimal argument $\mu_{k+1} \in \Lambda_{h}$ in step (ii) of Algorithm \ref{Alg:Majorant} is searched in the finite dimensional space $\Lambda_{h} \subset \Lambda$ of piecewise constant functions on $\mathcal{T}$. Then, under the assumption of $\phi \in V_h$, $f \in \Lambda_{h}$ with given values
$$ \phi(x_j), \phi(x_{j+1}), \quad v(x_j), v(x_{j+1}), \quad f(x_{j+\frac{1}{2}}), \quad (\tau^{*}_{k+1})'(x_{j+\frac{1}{2}}) $$
for $j=1\hdots n-1$, we obtain the formula
\begin{equation}
\mu_{k+1}(x_{j+\frac{1}{2}})=\left(-(\tau^{*}_{k+1})'(x_{j+\frac{1}{2}})-f(x_{j+\frac{1}{2}}) -\frac{v(x_{j})+v(x_{j+1})-\phi(x_{j})-\phi(x_{j+1})}
{2\left(1+\frac{1}{\beta_{k}}\right)} \right)^+,
\end{equation}
where $(\cdot)^+=\max\{ 0,\cdot \}$. Some iterations of Algorithm \ref{Alg:Majorant} are displayed in Figure \ref{Fig:Majorant}. We use a high (10000 in all experiments) number of iterations  in order to achieve the sharpest possible estimate \eqref{majorantEstimate}. Algorithm \ref{Alg:Majorant} provides a high quality approximations
$\tau^* \in V_h$ and $\lambda \in \Lambda_h$ in accordance with Remark \ref{rem:optimal}. We note that  Algorithm \ref{Alg:Majorant} provides a sharp estimate \eqref{majorantEstimate} for all types a), b), c) of approximations $v \in V_{0,h}$. It corresponds to values around $1.00$ in the last column of Tables \ref{tab:diferent_f_HundredIterations}, \ref{tab:diferent_f_ThousandIterations}, \ref{tab:diferent_f_tenThousandIterations}.

\begin{remark}[Update of $\beta$]
The experiments showed that the update of $\beta$ in the step (iii) of Algorithm \ref{Alg:Majorant} should not be called in every iteration. It turns out useful to call  steps (i) and (ii) repeatedly and run step (iii) only after variables $\tau^*$ and $\mu$ stabilize. We updated $\beta$ during the 5000th and the final 10000th iterations.
\end{remark}

\begin{table}
\begin{center}
\begin{small}
\begin{tabular}{|r|r|r|r|r|r|}
\hline
$f$  & $\frac{1}{2}\|v-u\|^{2}_{E}$ & $J(v)-J(u)$ & $\sqrt{\frac{J(v)-J(u)}{\frac{1}{2}\|v-u\|^{2}_{E}}}$ & $\mathcal{M}(v,\dots)$ & $\sqrt{\frac{\mathcal{M}(v,\dots)}{J(v)-J(u)}}$ \\
\hline
-5& 2.54e-006& 2.54e-006& 1.00& 2.55e-006& 1.00\\
-6& 3.66e-006& 3.66e-006& 1.00& 3.67e-006& 1.00\\
-7& 4.98e-006& 4.98e-006& 1.00& 4.99e-006& 1.00\\
-8& 6.51e-006& 6.51e-006& 1.00& 6.52e-006& 1.00\\
-9& 2.27e-005& 2.27e-005& 1.00& 2.39e-005& 1.03\\
-10& 6.49e-005& 6.86e-005& 1.03& 7.41e-005& 1.04\\
-11& 8.25e-005& 9.91e-005& 1.10& 1.06e-004& 1.04\\
-12& 8.57e-005& 9.99e-005& 1.08& 1.07e-004& 1.04\\
-13& 8.42e-005& 1.13e-004& 1.16& 1.20e-004& 1.03\\
-14& 8.73e-005& 3.69e-004& 2.06& 3.75e-004& 1.01\\
-15& 8.86e-005& 6.44e-004& 2.70& 6.51e-004& 1.00\\
-16& 1.02e-004& 9.91e-004& 3.11& 9.98e-004& 1.00\\
-17& 1.13e-004& 1.26e-003& 3.34& 1.27e-003& 1.00\\
-18& 1.24e-004& 1.54e-003& 3.53& 1.55e-003& 1.00\\
-19& 1.36e-004& 1.73e-003& 3.57& 1.74e-003& 1.00\\
-20& 1.56e-004& 1.98e-003& 3.55& 1.99e-003& 1.00\\
\hline
\end{tabular}
\caption{Verification of majorant and energy estimates for problems with various $f$ computed on an uniform mesh with $641$ nodes. Discrete solutions $v$ is computed by 100 iterations of the Uzawa algorithm.}
\label{tab:diferent_f_HundredIterations}
\vspace{1cm}
\begin{tabular}{|r|r|r|r|r|r|}
\hline
$f$  & $\frac{1}{2}\|v-u\|^{2}_{E}$ & $J(v)-J(u)$ & $\sqrt{\frac{J(v)-J(u)}{\frac{1}{2}\|v-u\|^{2}_{E}}}$ & $\mathcal{M}(v,\dots)$ & $\sqrt{\frac{\mathcal{M}(v,\dots)}{J(v)-J(u)}}$ \\
\hline
-5& 2.54e-006& 2.54e-006& 1.00& 2.55e-006& 1.00\\
-6& 3.66e-006& 3.66e-006& 1.00& 3.67e-006& 1.00\\
-7& 4.98e-006& 4.98e-006& 1.00& 4.99e-006& 1.00\\
-8& 6.51e-006& 6.51e-006& 1.00& 6.52e-006& 1.00\\
-9& 9.04e-006& 9.49e-006& 1.02& 9.83e-006& 1.02\\
-10& 1.00e-005& 2.36e-005& 1.53& 2.39e-005& 1.01\\
-11& 1.19e-005& 2.72e-005& 1.51& 2.76e-005& 1.01\\
-12& 1.33e-005& 2.90e-005& 1.48& 2.94e-005& 1.01\\
-13& 1.53e-005& 3.96e-005& 1.61& 4.01e-005& 1.01\\
-14& 1.73e-005& 4.79e-005& 1.66& 4.85e-005& 1.01\\
-15& 1.88e-005& 5.61e-005& 1.73& 5.69e-005& 1.01\\
-16& 2.07e-005& 5.51e-005& 1.63& 5.59e-005& 1.01\\
-17& 2.24e-005& 6.01e-005& 1.64& 6.10e-005& 1.01\\
-18& 2.48e-005& 6.94e-005& 1.67& 7.03e-005& 1.01\\
-19& 2.70e-005& 9.06e-005& 1.83& 9.17e-005& 1.01\\
-20& 2.91e-005& 8.52e-005& 1.71& 8.63e-005& 1.01\\
\hline
\end{tabular}
\caption{Verification of majorant and energy estimates for problems with various $f$ computed on an uniform mesh with $641$ nodes. Discrete solutions $v$ is computed by 1000 iterations of the Uzawa algorithm.}
\label{tab:diferent_f_ThousandIterations}
\end{small}
\end{center}
\end{table}
\begin{table}
\begin{center}
\begin{small}
\begin{tabular}{|r|r|r|r|r|r|}
\hline
$f$  & $\frac{1}{2}\|v-u\|^{2}_{E}$ & $J(v)-J(u)$ & $\sqrt{\frac{J(v)-J(u)}{\frac{1}{2}\|v-u\|^{2}_{E}}}$ & $\mathcal{M}(v,\dots)$ & $\sqrt{\frac{\mathcal{M}(v,\dots)}{J(v)-J(u)}}$ \\
\hline
-5& 2.54e-006& 2.54e-006& 1.00& 2.54e-006& 1.00\\
-6& 3.66e-006& 3.66e-006& 1.00& 3.66e-006& 1.00\\
-7& 4.98e-006& 4.98e-006& 1.00& 4.98e-006& 1.00\\
-8& 6.51e-006& 6.51e-006& 1.00& 6.51e-006& 1.00\\
-9& 7.78e-006& 8.05e-006& 1.02& 8.10e-006& 1.00\\
-10& 9.11e-006& 9.79e-006& 1.04& 9.88e-006& 1.00\\
-11& 1.05e-005& 1.10e-005& 1.02& 1.11e-005& 1.00\\
-12& 1.20e-005& 1.29e-005& 1.04& 1.30e-005& 1.00\\
-13& 1.35e-005& 1.42e-005& 1.03& 1.44e-005& 1.00\\
-14& 1.51e-005& 1.60e-005& 1.03& 1.61e-005& 1.00\\
-15& 1.68e-005& 1.78e-005& 1.03& 1.79e-005& 1.00\\
-16& 1.84e-005& 2.01e-005& 1.04& 2.03e-005& 1.01\\
-17& 2.02e-005& 2.16e-005& 1.03& 2.18e-005& 1.00\\
-18& 2.20e-005& 2.39e-005& 1.04& 2.42e-005& 1.01\\
-19& 2.39e-005& 2.55e-005& 1.03& 2.57e-005& 1.00\\
-20& 2.58e-005& 2.80e-005& 1.04& 2.83e-005& 1.01\\
\hline
\end{tabular}
\caption{Verification of majorant and energy estimates for problems with various $f$ computed on an uniform mesh with $641$ nodes. Discrete solutions $v$ is computed by 10000 iterations of the Uzawa algorithm. }
\label{tab:diferent_f_tenThousandIterations}
\end{small}
\end{center}
\end{table}

\newpage
\section*{Conclusions}
A new minimization majorant algorithm providing an optimal value of the functional majorant $\mathcal{M}$ that bounds the difference of energies $J(v)-J(u)$ was described.
Numerical experiments in 1D show that the bound can be computed sharply for both low and high quality approximation $v$ assuming a high number of the algorithm iterations.
An analysis of a nonlinear benchmark with known analytical solution indicates, that $J(v)-J(u)$ provides the exact value of the error of approximation $\frac{1}{2}\|v-u\|^{2}_{E}$ in situations when contact zone of the discrete solution $v$ covers whole contact zone of the exact solution $u$.

\section*{Acknowledgment}
Both authors acknowledge the support of the European Regional Development Fund in the Centre of Excellence project IT4Innovations (CZ.1.05/1.1.00/02.0070) and by the project SPOMECH - Creating a multidisciplinary R\&D team for reliable solution of mechanical problems, reg. no. CZ.1.07/2.3.00/20.0070 within Operational Programme 'Education for competitiveness'  funded by Structural Funds of the European Union and  state budget of the Czech Republic. Authors would also like to thank to S. Repin (St. Petersburg), J. Kraus (Linz), D. Pauly (Duisburg-Essen) and O. Vlach (Ostrava) for discussions. Number of discussions with S. Repin were held during the visit of the second author at the Steklov institute of mathematics in St. Peterburg in frames of scientific cooperations of Czech and Russian academies of sciences - Institute of Geonics AS CR and St. Petersburg Department of Steklov Institute of Mathematics.

\noindent P. Harasim and J. Valdman \\
Centre of Excellence IT4Innovations, \\
V\v SB-Technical University of Ostrava, \\
Czech Republic \\
e-mail: jan.valdman@vsb.cz \\

\noindent P. Harasim \\
Faculty of Civil Engineering  \\
Brno University of Technology \\
Czech Republic \\

\noindent J. Valdman \\
Institute of Information Theory and Automation of the ASCR \\
Prague \\
Czech Republic\\


\begin{thebibliography}{5}
\bibitem{AiOd}
{\it M.~Ainsworth and J. T.~Oden},
A posteriori error estimation in finite element analysis,
Wiley and Sons, New York, 2000.

\bibitem{BabSt}
{\it I.~Babu{\v s}ka and T.~Strouboulis} \,:
The finite element method and its reliability,
Oxford University Press, New York, 2001.

\bibitem{BanRa}
W. Bangerth and R. Rannacher,
{\em Adaptive finite element methods for differential equations},
Birkh\"auser, Berlin, 2003.

\bibitem{BHS} {\it D. Braess, R. H. W. Hoppe, J. Sch\"{o}berl}\,: A posteriori estimators for obstacle problems by the hypercircle method. Comp. Visual. Sci. 11, 2008, 351--362.

\bibitem{BHR}
{\it F. Brezi, W. W. Hager, P. A. Raviart}\,:
Error estimates for the finite element solution of variational inequalities I.
Numer. Math., 28, 1977, 431--443.

\bibitem{BuRe}
{\it H.~Buss and S.~Repin}\,:
A posteriori error estimates for boundary value problems with obstacles,
Proceedings of 3nd European Conference on Numerical Mathematics and Advanced Applications, J\"{y}vaskyl\"{a}, 1999, World Scientific, 2000, 162–-170.

\bibitem{CaMe}
{\it C. Carstensen and C. Merdon},
A posteriori error estimator completition for conforming obstacle problems,
Numer. Methods Partial Differential Eq. 29, 2013, 667-–692.

\bibitem{DO}
{\it Z.~Dost\'al}:
Optimal Quadratic Programming Algorithms.
Springer 2009.

\bibitem{FALK}
{\it R. S. Falk}\,:
Error estimates for the approximation of a class of variational inequalities.
Math. Comput., 28, 1974, 963--971.

\bibitem{FuRe}
{\it M. Fuchs and S. Repin}\,:
A Posteriori Error Estimates for the Approximations of the Stresses in the Hencky Plasticity Problem.
Numer. Funct. Analysis and Optimization, 32, 2011.

\bibitem{GLT}
{\it R. Glowinski, J. L. Lions, R. Tr\'emolieres}\,:
Numerical analysis of variational inequalities. North-Holland, 1981.

\bibitem{HHNL}
{\it I. Hlav\'a\v cek, J. Haslinger, J. Ne\v cas, J. Lov\'i\v sek}\,:
Solution of variational inequalities in mechanics, vol. 66. Applied Mathematical Sciences, Springer-Verlag, New York, 1988.

\bibitem{KO}
{\it N. Kikuchi, J. T. Oden}:
Contact Problems in Elasticity: A Study of Variational Inequalities and Finite Element Methods,
SIAM, 1995.

\bibitem{KrTo}
{\it J. Kraus, S. Tomar}\,: Algebraic multilevel iteration method for lowest-order Raviart-Thomas space and applications.
Int. J. Numer. Meth. Engng. 86, 2011, 1175--1196.

\bibitem{LiSt}
{\it J. L. Lions, G. Stampacchia}\,:
Variational inequalities.
Comm. Pure Appl. Math., XX(3), 1967, 493--519.

\bibitem{NeRep} {\it P. Neittaanm\"{a}ki, S. Repin}\,:
Reliable methods for computer simulation (error control and a posteriori estimates). Elsevier, 2004.

\bibitem{Re1}
{\it S.~Repin}\,:
A posteriori error estimation for variational problems with uniformly convex functionals,
{\em Math. Comput.} 69, No. 230, 2000, 481--500.

\bibitem{Re2}
{\it S.~Repin}\,:
A posteriori error estimation for nonlinear variational problems by duality theory,
Zapiski Nauchn. Semin. POMI 243, 1997, 201--214.

\bibitem{Re4}
{\it S.~Repin}\,:
Estimates of deviations from exact solutions of elliptic variational inequalities,
Zapiski Nauchn. Semin, POMI 271, 2000, 188--203.

\bibitem{ReGruyter}
{\it S.~Repin}\,:
A posteriori estimates for partial differential equations,
Walter de Gruyter, Berlin, 2008.

\bibitem{ReVa}
{\it S. Repin, J. Valdman}\,:
Functional a posteriori error estimates for problems with nonlinear boundary conditions,
Journal of Numerical Mathematics 16, No. 1, 2008, 51--81.

\bibitem{ReVa2}
{\it S. Repin, J. Valdman}\,:
Functional a posteriori error estimates for incremental models in elasto-plasticity.
Cent. Eur. J. Math. 7, No. 3, 2009, 506--519.

\bibitem{Ul}
{\it M. Ulbrich}\,:
Semismooth Newton methods for variational inequalities and constrained optimization problems in function spaces.
SIAM 2011.

\bibitem{Vld}
{\it J. Valdman}\,:
Minimization of functional majorant in a posteriori error analysis based on $H(div)$ multigrid-preconditioned CG method.
Advances in Numerical Analysis, 2009.

\bibitem{ZVKG}
{\it Q. Zou, A. Veeser, R. Kornhuber, C. Gr\"aser}\,:
Hierarchical error estimates for the energy functional in obstacle problems,
Numerische Mathematik 117, No. 4, 2012, 653--677.

\end{thebibliography}
\end{document}